\documentclass[11pt]{article}
\usepackage{latexsym, amsmath, amssymb, amsfonts, amscd}
\usepackage{graphics}
\usepackage[dvips]{epsfig}

\newtheorem{thm}{Theorem}[section]
\newtheorem{prop}[thm]{Proposition}
\newtheorem{lemma}[thm]{Lemma}

\newtheorem{dfn}[thm]{Definition}

\newtheorem{rmk}[thm]{Remark}

\newcommand{\reals}{\mathbb R}

\newcommand{\call}{{\cal L}}

\def\qed{\rule{2.3mm}{2.3mm}}

\newcommand{\half}{\frac{1}{2}}

\def\lcf{\lbrack\! \lbrack}
\def\rcf{\rbrack\! \rbrack}
\def\gpd{\rightrightarrows}
\def\gr{G \underset{\t} {\overset {\s} \gpd} M}

\def\s{\alpha}
\def\t{\beta}

\begin{document}
\title{\bf Integration of Dirac-Jacobi structures}
\author{
David Iglesias Ponte\thanks{e-mail: iglesias@math.psu.edu} \ and
A\"{\i}ssa Wade\thanks{e-mail: wade@math.psu.edu} \footnote{The
first author is partially supported by
MCYT grant BFM2003-01319} \\
 {\small Department of Mathematics, Penn State University} \\
}

\date{}
\maketitle
\begin{abstract}
We study precontact groupoids whose infinitesimal counterparts are
Dirac-Jacobi structures. These geometric objects generalize
 contact groupoids. We also explain the relationship
between precontact groupoids and  homogeneous presymplectic
groupoids. Finally, we present some examples of precontact
groupoids.
\end{abstract}

\section{Introduction}
Presymplectic groupoids, introduced and studied in \cite{BCWZ}, are
global counterparts of  Dirac structures. They  allow to extend the
well-known correspondence between symplectic groupoids and Poisson
manifolds to the context of Dirac geometry. Moreover, they provide a
framework for a unified formulation of various notions of momentum
maps (\cite{BC}). On the other hand, Dirac-Jacobi structures (called
${\cal E}^1(M)$-Dirac structures in \cite{Wa00})
 include both Dirac and Jacobi structures.
They naturally appeared in the geometric prequantization of Dirac
manifolds \cite{WZ04,ZZ05}.

\bigskip

In this paper, our aim is to investigate the integrability problem
for Dirac-Jacobi structures. This work is motivated
 by the fact that many Dirac manifolds can
be quantized through their integrating Lie groupoids. We show that
the global counterparts of  Dirac-Jacobi manifolds are what we call
here precontact groupoids. In particular, we recover the
integrability of Jacobi structures \cite{CZ04}. We also prove that
there is a one-to-one correspondence between precontact groupoids
and homogeneous presymplectic groupoids. Moreover,  the precontact
groupoid $\widetilde{G}$ associated with an integrable Dirac
structure $L_0$ on $M$ is just the prequantization of the
presymplectic groupoid $G$ associated with $L_0$ (that is, the
central extension of Lie groupoids $M\times S^1 \rightarrow
\widetilde{G} \rightarrow G$ satisfying some compatibility
conditions), provided that the canonical Dirac Jacobi
 structure $L$ on $M$  corresponding $L_0$ is integrable, see
Section \ref{Ex-Dirac}. We should mention that M. Zambon and C. Zhu
independently study the geometry of prequantization spaces
\cite{ZZ05}.

Here is an outline of the paper. In Sections 2 and 3, we give some
definitions and results needed to establish our results. Section 4
contains our main results (Theorems \ref{thm1} and \ref{thm2}). In
Section 5, we give  some examples of precontact groupoids.

\section {Basic definitions and results}
\subsection{Dirac structures and presymplectic groupoids}
Let $M$ be a smooth $n$-dimensional manifold. There is a natural
symmetric pairing $ \langle \cdot ,\cdot \rangle$ on the vector
bundle $TM\oplus T^*M$ given by
\[
\langle X_1 +\xi _1,X_2+\xi _2\rangle = \half \Big ( \xi _1
(X_2)+\xi _2 (X_1) \Big ).
\]
Furthermore, the space of smooth sections of $TM\oplus T^\ast M$ is
endowed with the Courant bracket, which is defined by
\[
[X_1 +\xi _1,X_2+\xi _2]=[X_1 ,X_2]+{\cal L}_{X_1} \xi _2-
i_{X_2}d\xi _1 ,
\]
\noindent for any $X_1 +\xi _1, \ X_2 +\xi _2 \in  \Gamma(TM\oplus
T^*M).$
\begin{dfn}{\rm \cite{BCWZ,C}} A \textbf{Dirac structure} on a smooth
manifold $M$ is a subbundle $L$ of $TM\oplus T^\ast M$ which is
maximally isotropic with respect to the symmetric pairing $ \langle
\cdot ,\cdot \rangle$ and whose space of sections is closed under
the Courant bracket.

\smallskip
Let $L_M$ and $L_N$ be Dirac structures on $M$ and $N$,
respectively. We say that a smooth map $F:M\to N$ is a (forward)
\textbf{Dirac map} if $L_N=F_\ast (L_M)$, where
\[
F_\ast(L_M)=\{ (dF)(X)+\xi \ | \ X+F^\ast \xi\in L_M \}.
\]
\end{dfn}
Recall that any Dirac structure $L$ has an induced Lie algebroid
structure: the Lie bracket on $\Gamma (L)$ is just the restriction
of the Courant bracket and the anchor map is the restriction  of the
first projection to $L$, i.e. $pr_1|_{L}: L \to TM$. Now, suppose
that $L$ is a Dirac structure which is isomorphic to the Lie
algebroid  of a Lie groupoid  $G$. Such a Lie groupoid is called an
{\bf integration} of the Lie algebroid $L$. Then, there exists an
induced closed 2-form  on $G$ with some
 additional properties. More precisely,
\begin{dfn}\label{presymplectic}{\rm \cite{BCWZ}}
A \textbf{presymplectic groupoid} is a pair $(G,\omega )$ which
consists of a groupoid $\gr$ such that $\mbox{\rm dim(G)}=2\,
\mbox{\rm dim(M)}$, and a 2-form $\omega\in \Omega ^2(G)$ satisfying
the following conditions:
\begin{itemize}
\item[{\rm (i)}] $\omega$ is closed, i.e. $d\omega =0$.
\item[{\rm (ii)}] $\omega$ is multiplicative, that is, $m^\ast
\omega =pr_1^\ast \omega +pr_2^\ast \omega$.
\item[{\rm (iii)}] $\mbox{Ker }(\omega _x)\cap
\mbox{Ker }(d\s)_x\cap\mbox{Ker }(d\t)_x=\{0\}$, for all $x\in M$.
\end{itemize}
We  say that the presymplectic groupoid $(G,\omega )$ is
\textbf{homogeneous} if there exists a mutiplicative vector field
$Z$ such that ${\cal L}_Z \omega =\omega$.
\end{dfn}
The relationship between Dirac structures and presymplectic
groupoids is provided by the following result:
\begin{prop}\label{Prop 1}{\rm \cite{BCWZ}}
Given  a presymplectic groupoid $(G ,\omega )$, there is a canonical
Dirac structure $L$ on $M$ which is isomorphic to the Lie algebroid
$AG$ of $G$, and such that the target map $\t :(G, L_\omega) \to (M,
L)$ is a Dirac map, while the source map $\s: (G, L_\omega) \to (M,
L)$ is anti-Dirac.

Conversely, suppose that $L$ is a Dirac structure on $M$ whose
associated Lie algebroid is integrable, and let $G(L)$ be its
$\s$-simply connected integration. Then, there exists a unique
2-form $\omega$ such that $(G(L),\omega )$ is a presymplectic
groupoid, the target map is a Dirac map, and the source map is
anti-Dirac.
\end{prop}

\subsection{Dirac-Jacobi structures}
Let $M$ be a smooth $n$-dimensional manifold. There is a natural
bilinear operation $ \langle \cdot ,\cdot \rangle$ on the vector
bundle ${\cal E}^1(M)= (TM \times \reals) \oplus (T^*M \times
\reals)$ defined by:
$$\Big \langle(X_1, f_1)+ (\xi_1, g_1), (X_2, f_2)+ (\xi_2, g_2)
 \Big \rangle =\frac{1}{2} (i_{X_2} \xi_1 +
i_{X_1} \xi_2 +f_1 g_2 + f_2 g_1),$$ for any $(X_{\ell}, f_{\ell})+
(\xi_{\ell}, g_{\ell}) \in \Gamma({\cal E}^1(M))$, with $\ell=1,2$.
In addition, the space of smooth sections of ${\cal E}^1(M)$ is
equipped with an $\reals$-bilinear operation which can be viewed as
an extension of the Courant bracket on $TM \oplus T^*M$, i.e.
\[
\begin{array}{ccl}
\makebox{{\bf [}}(X_1,f_1)+(\xi _1,g_1),(X_2,f_2)+(\xi
_2,g_2)\makebox{{\bf ]}}&=&\Big ( [X_1,X_2],X_1(f_2)-X_2(f_1)\Big
)\\ &&\kern-110pt+\Big ({\cal L}_{X_1}\xi _2-{\cal L}_{X_2}\xi
_1+\frac{1}{2}d(i_{X_2}\xi _1-i_{X_1}\xi _2)\\
&&\kern-115pt+f_1\xi _2-f_2\xi
_1+\frac{1}{2}(g_2df_1-g_1df_2-f_1dg_2+f_2dg_1),\\
&&\kern-120ptX_1(g_2)-X_2(g_1)+\frac{1}{2}(i_{X_2}\xi _1-i_{X_1}\xi
_2-f_2g_1+f_1g_2)\Big ),
\end{array}
\]
for any  $(X_{\ell},f_{\ell})+(\xi_{\ell},g_{\ell}) \in \Gamma({\cal
E}^1(M))$ with $\ell=1,2$. For an alternative description of this
bracket, see \cite{IW05}.
\begin{dfn}{\rm \cite{Wa00}}
A \textbf{Dirac-Jacobi structure} is a subbundle $L$ of ${\cal
E}^1(M)$ which is maximally isotropic with respect to $ \langle
\cdot,\cdot  \rangle$  and such that $\Gamma(L)$ is closed under the
extended Courant bracket $\makebox{{\bf [}}\cdot , \cdot
\makebox{{\bf ]}}$.

\smallskip
Let $L_M$ (resp., $L_N$) be a Dirac-Jacobi structure on $M$ (resp.,
$N$). We say that a smooth surjective map $F:M\to N$ is a (forward)
\textbf{Dirac-Jacobi map} if $L_N=F_\ast (L_M)$, where
\[
F_\ast(L_M)=\{((dF)(X),f)+(\xi, g) \ | \ (X, f\circ F)+(F^\ast \xi,
g \circ F)\in L_M\}.
\]
\end{dfn}
Basic examples of Dirac-Jacobi structures are Dirac and Jacobi
structures on $M$ (this explains the terminology introduced in
\cite{GM03}).
\subsection{Action Lie algebroids and 1-cocycles}
It is known that, given any a Lie algebroid $(A, \lcf \cdot , \cdot
\rcf, \rho)$ over $M$ and any 1-cocycle $\phi \in \Gamma(A^\ast)$,
there is an associated Lie algebroid  over $M \times \reals$,
denoted by $(A\times_{\phi} \reals, \lcf \cdot , \cdot \rcf ^\phi
,\rho ^\phi)$, where the smooth sections of $A\times_{\phi}\reals$
are of the form $\bar{X}(x,t)= X_t(x)$, with $X_t \in \Gamma (A)$
for all $t\in \reals$, and
\begin{equation}\label{corchbarra}
\begin{array}{l}
\lcf \bar{X},\bar{Y}\rcf  ^\phi(x,t) = \lcf X_t, Y_t \rcf (x)+ \phi
(X_t)(x) \frac{\partial \bar{Y}}{\partial t}-\phi
(Y_t)(x)\frac{\partial
\bar{X}}{\partial t},\\[5pt]
\rho^\phi (\bar{X} )(x,t)= \rho (X_t)(x) +\phi (X_t)(x)
\frac{\partial}{\partial t},
\end{array}
\end{equation}
where $\displaystyle \frac{\partial\bar{X}}{\partial t}\in \Gamma
(A\times_{\phi} \reals)$ denotes the derivative of $\bar{X}$ with
respect to $t$.

\begin{rmk}
\label{note} {\rm If $L$ is a Dirac-Jacobi structure then the
restriction of the extended Courant bracket to sections of $L$
together with the canonical projection of $L$ onto $TM$ make $L$
into a Lie algebroid over $M$. In addition, $\phi \in
\Gamma(L^\ast)$ defined by
\begin{equation}\label{1-cociclo}
\phi (v)=f,\mbox{ for } v=(X,f)+(\xi ,g)\in \Gamma (L),
\end{equation}
is a 1-cocycle for the Lie algebroid cohomology (see \cite{IM02}).
On the other hand, it is known that a Dirac-Jacobi structure $L$ on
$M$ corresponds to a Dirac structure $\widetilde{L}$ on $M \times
\reals$ given by
\begin{equation}\label{induced-Dirac}
\widetilde{L}=\Big \{\Big ( X+f \frac{\partial}{\partial t}\Big
)+\Big ( e^t(\xi+ g\,dt)\Big ) \,\Big | \, (X, f)+(\xi, g)\in L\Big
\}.
\end{equation}
Moreover, $\widetilde{L}$ is isomorphic to $L\times _\phi \reals$.}
\end{rmk}

\subsection{Conformal classes of Dirac-Jacobi structures}
 Let $L$ be a  Dirac-Jacobi structure on $M$ and let
$\varphi$ be a smooth nowhere vanishing function on $M$. We set
$\mu= d \ln|\varphi|.$ Consider the vector bundle $L_{\varphi}$ over
$M$ whose space of smooth sections is given by
$$\Gamma (L_{\varphi})= \{ (X, f-\mu(X))+\varphi (\xi + g \mu, g) \ | \
 (X,f) +(\xi, g) \in \Gamma(L) \}.$$
One can easily check that $L_{\varphi}$ is also a Dirac-Jacobi
structure on $M$. The correspondence $(L, \varphi) \mapsto
L_{\varphi}$ is called a conformal change. For any fixed $L$, the
family of all $(L_{\varphi})$ is called a conformal class of
Dirac-Jacobi structures.  For instance, when $L$ comes from a
presymplectic  form $\omega$ then $ L_{\varphi}$ is nothing but the
Dirac-Jacobi structure associated with $(\varphi \omega, d \ln |
\varphi |)$ (see \cite{Wa00, Wa02} for more details).

\section{Precontact groupoids}
\begin{dfn}
Let $\gr$ be a Lie groupoid such that $\mbox{\rm dim(G)}=2\,
\mbox{\rm dim(M)}+1$. A \textbf{precontact groupoid structure} on
$G$ is given by a pair $(\eta, \sigma)$ consisting of a 1-form
$\eta$ and a multiplicative function $\sigma$ (i.e., $\sigma
(gh)=\sigma (g)+\sigma (h)$) such that
\begin{equation}\label{sigma-mult}
m^*\eta= pr_1^* \eta+pr^*_1(e^{\sigma} ) pr_2^* \eta .
\end{equation}
\begin{equation}\label{non-deg}
\mbox{Ker }(d\eta _x)\cap \mbox{Ker }(\eta _x)\cap
 \mbox{Ker }(d\s)_x\cap\mbox{Ker }(d\t)_x=\{0\},
\mbox{ for all }x\in M.
\end{equation}
Two precontact structures $(\eta, \sigma)$ and $(\eta', \sigma')$ on
$G$ are equivalent if  there exists a nowhere vanishing function
$\varphi : M \to \reals$ such that
\[
\eta'= (\varphi \circ \s) \ \eta, \quad \quad
 \sigma'= \sigma+ \ln \Big| \frac{\varphi \circ \s}{\varphi \circ \t}\Big|
 \]
\end{dfn}
Now, consider a Lie groupoid  $\gr$ together with  a multiplicative
function $\sigma$. One can define a right action of $G$ on the
canonical projection $\pi _1:M\times \reals  \to M$ as follows:
\[
(x,t)\cdot g=(\s (g),t+\sigma (g)),\mbox{ for } (x,t,g)\in M\times
\reals \times G, \mbox{ such that }\t(g)=x.
\]
Therefore, we have the corresponding action groupoid $G\times
\reals\gpd M\times \reals$, denoted by $G\times _\sigma \reals $,
with structural functions given by
\begin{equation}\label{groupoid-GxR}
\begin{array}{l}
\s_\sigma(g,t)=(\s(g),\sigma (g)+ t), \\
\t_\sigma(h,s)=(\t(h),
s),\\
m_\sigma((g,t),(h,s))=(gh,t),\mbox{ if } \s _\sigma(g,t)=\t
_\sigma(h,s).
\end{array}
\end{equation}
We denote by $(AG,\lcf\;,\;\rcf,\rho)$  the Lie algebroid of $G$.
The multiplicative function $\sigma$ induces a $1$-cocycle $\phi$ on
$AG$ given by
\begin{equation}\label{induced-cocycle}
\phi(x)(X_x)=X_x(\sigma),\;\;\; \mbox{for $x\in M$ and $X_x\in A_x
G.$}
\end{equation}
In addition, we can identify the Lie algebroid of the Lie groupoid
$G\times_\sigma \reals$ with $AG\times _\phi \reals$.
 Conversely, one has the following
\begin{prop}\label{gpd-relation}{\rm \cite{CZ04}}
Let $L$ be a Lie algebroid over $M$, $\phi$ be a 1-cocycle and
$L\times _\phi \reals$ the Lie algebroid given by Equation
(\ref{corchbarra}). Then, $L$ is integrable if and only $L\times
_\phi \reals$ is integrable. Moreover, if $G(L)$ (resp., $G(L\times
_\phi \reals)$) is the $\s$-simply connected integration of $L$
(resp., $L\times _\phi \reals$) and $\sigma$ is the multiplicative
function associated with $\phi$, then $G(L\times _\phi \reals) \cong
G(L)\times_\sigma \reals$.
\end{prop}
There is a correspondence between pre\-contact and pre\-symplec\-tic
group\-oids. Indeed, one has the following proposition:
\begin{prop}\label{prep1}
Let $G  \gpd M$ be a Lie groupoid and $\sigma$ a multiplicative
function on $G$. There is a one-to-one correspondence between
precontact groupoids on $(G, \sigma)$ and homogeneous presymplectic
groupoids on $G \times_{\sigma} \reals$.
\end{prop}
{\it Proof:} We know that there exists a one-to-one correspondence
between 1-forms on a manifold $M$ and 2-forms on $M\times \reals$
homogeneous with respect to $\frac{\partial}{\partial t}$. More
precisely, if $\eta$ is a 1-form on $M$ then $\omega= d(e^t \eta)$
is a homogeneous 2-form on $M\times \reals$. Conversely, assume that
$\omega$ is homogeneous and set $\eta=i_{\frac{\partial} {\partial
t}} \omega$. One can check that $\call_{\frac{\partial} {\partial
t}}\eta = \eta$. Hence $\eta$ can be identified with a 1-form on
$M$. Using the relation $\omega =d(e^t \eta)$, it is straightforward
to prove that conditions (ii) and (iii) in Definition
\ref{presymplectic} are equivalent to Equations (\ref{sigma-mult})
and (\ref{non-deg}). \hfill \qed

\section{Integration of Dirac-Jacobi structures}
In this section, we show that precontact groupoids are the global
objects corresponding to Dirac-Jacobi structures. First,  note that
the following lemma which is an immediate consequence of Remark
\ref{note} and Proposition \ref{gpd-relation}.
\begin{lemma}
A  Dirac-Jacobi structure  $L$ is integrable if and only if its
associated Dirac structure $\widetilde{L} \subset T(M \times
\reals)\oplus T^*(M\times \reals)$ is integrable.
\end{lemma}
\begin{thm}\label{thm1}
Let $L$ be an integrable Dirac-Jacobi structure on $M$, and let
$G(L)$ be its $\s$-simply connected integration. There exists  a
multiplicative function $\sigma$  and a 1-form $\eta$ such that
$(G(L), \eta, \sigma)$  is a precontact groupoid. Furthermore, any
conformal class of integrable Dirac-Jacobi structures on $M$ induces
a conformal class of precontact  groupoid structures.
\end{thm}
{\it Proof:} Suppose $L$ is an integrable Dirac-Jacobi structure on
$M$. We denote by  $\phi \in \Gamma(L^*)$ the 1-cocycle defined by
Equation (\ref{1-cociclo}). By integration, the 1-cocycle $\phi$
induces a multiplicative function $\sigma$ on $G(L)$. Denote by
$\widetilde{L}$ the Dirac structure on $M \times \reals$ associated
with $L$ given by Equation (\ref{induced-Dirac}). Applying
Proposition \ref{Prop 1}, one gets a presymplectic groupoid
$(G(\widetilde{L}), \omega)$. Moreover, one gets from Remark
\ref{note}
 and Proposition \ref{gpd-relation} that
\[
G(\widetilde{L}) \cong G(L\times_{\phi}\reals) \cong
G(L)\times_{\sigma} \reals.
\]
Observe that $\omega$ is homogeneous with respect to the
multiplicative vector field ${\frac{\partial} {\partial t}}$.
Indeed, the differentiable family of diffeomorphisms $\psi_s : M
\times \reals \to M \times \reals$ defined by $\psi_s (x,t)= (x,
s+t)$ consists of Dirac maps, i.e $(\psi_s)_* (\widetilde{L})=
\widetilde{L}.$ Furthermore, by integration and using the flow
$\psi_s : G(L) \times \reals \to G(L) \times \reals$ $(g,t)\mapsto
(g,t+s)$ of the vector field $\frac{\partial}{\partial t}$, we get
\[
\call_{\frac{\partial} {\partial t}} \omega = \omega.
\]
Thus, using Proposition \ref{prep1}, we deduce that there exists a
1-form $\eta$ such that $(G(L), \eta, \sigma)$ is precontact. Recall
that, for every nowhere vanishing function $\varphi$ on $M$, there
is an equivalent Dirac-Jacobi structure on $M$ whose space of
sections is given by
$$\Gamma (L_{\varphi})= \{ (X, f-\mu(X))+\varphi (\xi + g \mu, g) \ | \
 (X,f) +(\xi, g) \in \Gamma(L) \}.$$
Consider the 1-cocycle $\phi_{\varphi} \in \Gamma(L^*_{\varphi})$
defined as follows
\[
\phi_{\varphi}(e)= f- \mu(X), \quad \quad \mathrm{for \ all} \ e=(X,
f-\mu(X))+\varphi (\xi + g \mu, g).
\]
We have a natural commutative diagram of vector bundle morphisms:
\begin{displaymath}
\begin{matrix}
L & \stackrel{\simeq}  \longrightarrow  & L_{\varphi}  & & & & &\cr
&&  & & &  & &\cr \downarrow  & & \downarrow  & & & & &\cr && & & &
& & \cr L \times_{\phi} \reals & \stackrel \simeq  \longrightarrow
 & L_{\varphi} \times_{\phi_{\varphi}} \reals & & & & & \cr
\end{matrix}
\end{displaymath}
 Moreover, $L_{\varphi}$ induces a precontact  structure
 $(\eta_{\varphi}, \sigma_{\varphi})$  on $G(L_{\varphi}) \cong G(L)$
 given by
$$\eta_{\varphi}= (\varphi \circ \alpha )\ \eta \quad \quad
 \sigma_{\varphi}=\sigma + \ln \Big| \frac{\varphi \circ \s}{\varphi \circ \t}\Big|.$$
Thus, any conformal class of integrable Dirac-Jacobi structures on
$M$ induces a conformal class of precontact  groupoid structures.
\hfill \qed

\medskip
 Conversely, we have the following result.
\begin{thm}\label{thm2}
Let $(G, \eta, \sigma)$ is a precontact groupoid over $M$. Then,
there  exists a canonical Dirac-Jacobi structure $L_M$ on $M$ which
is isomorphic to the Lie algebroid $AG$. Moreover,
 $\t :G\to M$ is a Dirac-Jacobi map.
\end{thm}
{\it Proof:} By Proposition \ref{prep1}, one has
 the presymplectic groupoid
$(G\times _\sigma \reals, \omega = d(e^t\eta ))$ over $M\times
\reals$. Then, using Proposition \ref{Prop 1}, one gets a Dirac
structure $L_{M\times \reals}$ on $M\times \reals$ such that $\t
_\sigma$ is a Dirac map. Thus,
\[
\begin{array}{rcl}
L_{M\times \reals}&=&\Big \{ d\t _\sigma\Big(X+F {\frac{\partial}
{\partial t}}\Big ) + \Big ( \xi + G dt \Big ) \, \Big | \, i_{(X+ F
\frac{\partial} {\partial t})}\omega = \t _\sigma ^\ast (\xi + G dt
)\Big \}
\\[8pt] &=&\Big \{d\t _\sigma\Big(X+F
\frac{\partial} {\partial t}\Big ) + \Big ( \xi +G dt \Big ) \, \Big
| \, \t _\sigma^\ast (\xi )=e^t( i_{X}d\eta +F \eta ),\,\,
G=-e^t\eta (X) \Big \} .
\end{array}
\]
Therefore, one can write $L_{M\times \reals}\equiv \widetilde{L}_M$,
where $L_M$  is  the Dirac-Jacobi structure given by
\[
L_{M}=\Big \{ ( (d\t)(X), F  ) +  ( \xi , G ) \, \Big | \, \t ^\ast
(\xi )= i_{X}d\eta +F \,\eta ,\,\,  G=-\eta (X) \Big \} .
\]
By construction, $\t :G\to M$ is a Dirac-Jacobi map. There follows
the result. \hfill \qed

\section{Examples}
In this section, we will give some examples of Dirac-Jacobi
structures and describe their corresponding precontact groupoids.
\subsection{Precontact structures}
A \textbf{precontact structure} on a manifold $M$ is just a 1-form
$\theta$ on $M$. A precontact structure $\theta$ on $M$ induces a
Dirac-Jacobi structure $L_\theta$ whose  space of smooth sections is
\[
\Gamma (L_\theta)=\{ (X,f )+(i_Xd\theta +f\,\theta ,-i_X\theta ) \ |
\ (X,f )\in \mathfrak X(M)\times C^\infty (M) \}.
\]
We observe that the Lie algebroids $L_\theta$ and $(TM\times \reals
,\makebox{{\bf [}}\cdot ,\cdot \makebox{{\bf ]}}, \pi)$ are
isomorphic, where $\pi :TM\times \reals \to TM$ is the canonical
projection over the first factor and $\makebox{{\bf [}}\cdot ,\cdot
\makebox{{\bf ]}}$ is given by
\[
\makebox{{\bf [}} (X,f),(Y,g)\makebox{{\bf ]}}=([X,Y],X(g)-Y(f)),
\]
for $(X,f),(Y,g )\in \mathfrak X (M)\times C^\infty (M)$. Moreover,
under the isomorphism between $L_\theta$ and $TM\times\reals$, the
1-cocycle $\phi $ is the pair $(0,1)\in \Omega^1(M)\times C^\infty
(M)$.

On the other hand, consider the product $G=M\times M\times \reals$
of the pair groupoid with $\reals$. The function $\sigma :G\to
\reals$, $(x,y,t)\mapsto t$, is trivially multiplicative. In
addition, if  $\theta$ is a 1-form on $M$ then one can define the
1-form $\eta$ on $G$ given by
\[
\eta =\pi _1^\ast \theta - e^{\sigma }\pi _2^\ast \theta,
\]
where $\pi _i$, $i\in \{1,2\}$, is the projection on the $i$-th
component. Then $(G,\eta, \sigma)$ is a precontact groupoid and,
moreover, the corresponding Dirac-Jacobi structure on $M$ is just
$L_\theta$.

\subsection{Dirac structures}\label{Ex-Dirac}
Let $L_0$ be a vector subbundle of $TM\oplus T^\ast M$ and consider
the vector subbundle $L$ of ${\cal E}^1(M)$ whose sections are
\[
\Gamma (L)=\{ (X,0)+(\alpha ,f )\,|\, X+\alpha \in \Gamma (L_0) ,f
\in C^\infty (M) \}.
\]
Then, $L_0$ is a Dirac structure on $M$ if and only if $L$ is a
Dirac-Jacobi structure.

If we  denote by $L_0$ (resp., $L$) the Lie algebroid associated
with the Dirac structure (resp., the Dirac-Jacobi structure), then
we have that $L\equiv L_0\times \reals$. Moreover, a direct
computation shows that the bracket on $\Gamma (L)$ is given by
\[
\lcf ({\cal X}_1,f_1),({\cal X}_2,f_2)\rcf _L=\Big ( \lcf {\cal
X}_1,{\cal X}_2\rcf _{L_0},\rho _{L_0}({\cal X}_1)(f_2)-\rho_{L_0}
({\cal X}_2)(f_1)+\Omega _{L_0}({\cal X}_1,{\cal X}_2)\Big ),
\]
for $({\cal X}_1,f_1),({\cal X}_2,f_2)\in \Gamma (L)$, and where
$(\lcf \, ,\, \rcf _{L_0},\rho _{L_0})$ (resp., $(\lcf \, ,\, \rcf
_L,\rho _L)$) denotes the Lie algebroid structure on $L_0$ (resp.,
$L$) and $\Omega _{L_0}\in \Gamma (\wedge ^2L_0^\ast )$ is the
closed 2-section given by
\[
\Omega _{L_0}({\cal X}_1,{\cal X}_2)=\half \Big ( \xi _1 (X_2)-\xi
_2 (X_1) \Big ), \mbox{ for }{\cal X}_i=X_i+\xi _i \in \Gamma (L_0).
\]
Therefore, we have that the Lie algebroid structure on $L$ is just
the central extension of the Lie algebroid $L_0$ by the closed
2-section $\Omega _{L _0}$. On the other hand, from Equation
(\ref{1-cociclo}), one sees that the 1-cocycle $\phi$ identically
vanishes. If $L_0$ as well as $L$ are integrable and $(G,\omega )$
is the presymplectic groupoid associated with $L_0$ then one obtains
a \textbf{prequantization} of $(G,\omega)$, that is, a central
extension of Lie groupoids
\[
M\times S^1 \rightarrow \widetilde{G} \rightarrow G,
\]
and a multiplicative 1-form $\eta\in \Omega ^1(\widetilde{G})$
($m^\ast \eta =pr_1^\ast \eta +pr_2^\ast \eta$) which is a
connection 1-form for the principal $S^1$-bundle $\pi
:\widetilde{G}\to G$ and which satisfies $d\eta=\pi ^\ast \omega$
(see \cite{C04}).
\subsection{Jacobi manifolds}
Let $(\Lambda ,E)$ be Jacobi manifold and $L_{(\Lambda ,E)}$ the
corresponding Dirac-Jacobi structure
\[
L_{(\Lambda ,E)}=\{ (\Lambda ^\sharp (\alpha _x)+\lambda \,
E_x,-\alpha _x (E_x) )+(\alpha _x,\lambda )\,|\, (\alpha _x,\lambda
) \in T^\ast _xM\times \reals, \,x\in M \}.
\]
In this case,  the corresponding Dirac structure on $M\times \reals$
is the one coming from the Poissonization of the Jacobi structure
$(\Lambda ,E)$, i.e. $\Pi =e^{-t}(\Lambda +\frac{\partial}{\partial
t}\wedge E)$. Thus, the presymplectic groupoid is a honest
symplectic groupoid (see \cite{BCWZ}), and therefore, the precontact
groupoid structure integrating $L_{(\Lambda, E)}$ is a
\textbf{contact groupoid}, i.e., the 1-form defines a contact
structure on $G$. This result was first proved in \cite{D97} (see
also, \cite{CZ04,KS93}).
%\section{Special case: regular Dirac-Jacobi structures}
\bigskip
{\bf Acknowledgments:} D.~Iglesias wishes to thank the Spanish
Ministry of Education and Culture and Fulbright program for a
MECD/Fulbright postdoctoral grant.  A. Wade would like to thank the
MSRI for hospitality while this paper was being prepared.

\end{document}